\documentclass{article}
\usepackage{33tree}

\begin{document}
\title{On open flat sets in  spaces with bipolar comparison}
\author{Nina Lebedeva}

\nofootnote{Partially supported by RFBR grant 17-01-00128.}

\newcommand{\Addresses}{{\bigskip\footnotesize

  Nina Lebedeva, \par\nopagebreak\textsc{Math. Dept.
St. Petersburg State University,
Universitetsky pr., 28, 
Stary Peterhof, 
198504, Russia.}
  \par\nopagebreak
  \textsc{Steklov Institute,
27 Fontanka, St. Petersburg, 
191023, Russia.}
  \par\nopagebreak
  \textit{Email}: \texttt{lebed@pdmi.ras.ru}
}}

\date{}

\maketitle

\begin{abstract}
We show that if a Riemannian manifold satisfies (3,3)-bipolar comparisons and has an open flat subset then it is flat.
The same holds for  a version of MTW where the perpendicularity is dropped. 

In particular we get that the (3,3)-bipolar  comparison is strictly stronger than the Alexandrov comparison.
\end{abstract}

\section{Introduction}

We say that a metric space $X$ satisfies the \emph{$(k,l)$-bipolar comparison} if 
for any $a_0,a_1,\dots,a_k; b_0,b_1,\dots, b_l \in X$
there are points $\hat a_0, \hat a_1,\dots, \hat a_n, \hat b_0, \hat b_1,\dots,\hat b_n $ in the Hilbert space $\HH$ 
such that 
\[|\hat a_0-\hat b_0|_\HH=|a_0-b_0|_X,\quad
|\hat a_i-\hat a_0|_\HH=|a_i-a_0|_X,\quad
|\hat b_i-\hat b_0|_\HH=|b_i-b_0|_X\]
for any $i,j$ and
\[|\hat x-\hat y|_\HH\ge|x-y|_X\]
for any $x,y\in\{a_0, a_1,\dots, a_k, b_0, b_1,\dots, b_l\}$.

This  definition was introduced in \cite{LPZ}.
The class of compact length metric spaces satisfying $(k,0)$-bipolar comparison with $k\ge 2$ coincide with the class of Alexandrov spaces
with nonnegative curvature, (for $k=2$ it is just one of
the equivalent definitions, for arbitrary $k$ see \cite{AKP},
\cite{LS}).
In general
 $(k,l)$-bipolar comparisons
(with $k$ or $l\ge 2$) 
  for length metric spaces are  stronger conditions
 than nonnegative curvature condition and they describe some
 new interesting classes of spaces.
In
 particular, we prove in \cite{LPZ} that
for Riemannian manifolds
 $(4,1)$-bipolar comparison 
  is equivalent to
the conditions related to the continuity of optimal transport.
Also in \cite {LPZ} we together with coauthors describe classes of Riemannian 
manifolds 
satisfying $(k,l)$--bipolar comparisons for almost all
$k, l$ excepting $(2,3)$ and $(3,3)$-bipolar  comparisons.
In particular it was not known if $(3,3)$-bipolar comparison
differs from Alexandrov's comparison.
In this note the affirmative  answer is obtained as a corollary
of some rigidity result for spaces with $(3,3)$-bipolar comparison.
To formulate exact statements we
need some definitions and notations.

Let $M$ be a Riemannian manifold and $p\in M$.
The subset of tangent vectors $v\in\T_p$ such that there is a minimizing geodesic $[p\,q]$ in the direction of $v$ with length $|v|$ will be denoted as $\overline{\TIL}_p$.
The interior of $\overline{\TIL}_p$ is denoted by $\TIL_p$; it is called \emph{tangent injectivity locus} at $p$.
If at $\TIL_p$ is convex for any $p\in M$, then $M$ is called CTIL.

Riemannian manifold $M$ satisfies MTW if the following holds.
 For any point $p\in M$, any $W\in \TIL_p$ and  tangent vectors 
$X,Y\in \T_p$, such that $X\perp Y$ we have
\[\frac{\partial^4}{\partial^2s\,\partial^2t}\left|\exp_p(s\cdot X)-\exp_p(W+t\cdot Y)\right|_M^2\le 0\eqlbl{eq:MTW}\]
at $t=s=0$. 

This definition was introduced by
Xi-Nan Ma, Neil Trudinger and Xu-Jia Wang 
in \cite{MTW},
 Cedric Villani studied a synthetic version of this definition
 (\cite{MTW+CTIL}).
If the same inequality holds without the assumption $X\perp Y$
Riemannian manifold $M$ satisfies
MTW$^{\not\perp}$ \cite{FRV-Nec+Suf}.

MTW and CTIL are necessary condition for TCP
(transport continuity property).
In \cite{FRV-Nec+Suf}, Alessio Figalli, Ludovic Rifford and C\'edric Villani showed that
a strict version of CTIL and MTW provide a sufficient condition for TCP.
A compact Riemannian manifold $M$ is called TCP 
if for any two regular measures with density functions bounded away from zero and infinity the generalized solution of Monge--Amp\`{e}re equation provided by optimal transport 
is a genuine (continuous) solution.

Let us denote by 
$\mathcal M_{(k,l)}$ the class of smooth complete Riemannian
manifolds satisfying $(k,l)$--bipolar comparison
and by $\mathcal M_{\ge 0}$ the class of complete  Riemannian
manifolds with nonnegative sectional curvature.

%Now we are about to formulate results from\cite{LPZ} describing classes of manifolds satisfying bipolar 
%comparisons.% and formulate some new results.

It was mentioned above, that 
 $$\mathcal M_{\ge 0}=\mathcal M_{(k,0)}$$
 for $k\ge 2$ and 
 it is obvious from definition, that
 $$\mathcal M_{(k',l')}\subset \mathcal M_{(k,l)}$$
 if $k'\ge k$ and $l'\ge l$.
It is proven in \cite{LPZ}
that
$$\mathcal M_{\ge 0}=\mathcal M_{(2,2)}=\mathcal M_{(3,1)}$$
 and
 $$\mathcal M_{(4,1)}= \mathcal M_{(k,l)}$$
 for $k\ge 4$ and $l\ge 1$.
 The most interesting fact proven in \cite{LPZ}
 is that 
 $$\mathcal M_{(4,1)}= \mathcal M_{CTIL} \cap
  \mathcal M_{MTW^{\not\perp}},$$
  where
   $\mathcal M_{CTIL}$,
 $ \mathcal M_{MTW^{\not\perp}}$
  are classes of smooth Riemannian
manifolds satisfying 
CTIL and MTW$^{\not\perp}$ correspondingly.
In particular this implies that $\mathcal M_{(4,1)}\neq \mathcal M_{\ge 0}$.

In this paper we prove the following two results.
 
\begin{thm}{Theorem}\label{thm-33}
Let $M$ be a complete Riemannian manifold that satisfies (3,3)-bipolar comparison and contains a nonempty open flat subset. Then $M$ is flat.  
\end{thm}

\begin{thm}{Theorem}\label{thm-mtw}
Let $M$ be a complete Riemannian manifold that satisfies MTW$^{\not\perp}$
and contains a nonempty open flat subset. Then $M$ is flat.
\end{thm}

\begin{thm}{Corollary}
We have that $\mathcal M_{(3,3)}\neq \mathcal M_{\ge 0}$.
\end{thm}

Theorem~\ref{thm-33} follows from Proposition~\ref{prop-33}
and
Theorem~\ref{thm-mtw} follows from Proposition~\ref{prop-mtw}, proved in the next section.

As a related result we would like to mention
a rigidity result for manifolds with nonnegative sectional
curvature with flat open subsets
by Dmitrii Panov and Anton Petrunin \cite{PaPe}.

\section{Proofs}
For points $a,b,c$ in a manifold
we denote by
$\measuredangle[a\,{}^{b}_{c}]$ the angle at $a$ of the triangle
$[abc]$.

\begin{thm}{Key lemma}\label{lem:key}
Let $M$ be a complete Riemannian manifold that satisfies 
(3,3)-bipolar comparison.
Assume that for the points $x_p,p,q,x_q$ in $M$ there is a triangle $[\~p\~q\~x]$ in the 
Euclidean
plane $\EE^2$  such that 
\[|x_p-p|_M=|\~x-\~p|_{\EE^2},\quad |p-q|_M=|\~p-\~q|_{\EE^2},\quad |q-x_q|_M=|\~q-\~x|_{\EE^2}\] 
and moreover a neighborhood $N\subset\EE^2$ of the base $[\~p\~q]$ admits a globally isometric embedding $\iota$ into $M$ such that $\iota([\~p\~x]\cap N)\subset [px_p]$ and $\iota([\~q\~x]\cap N)\subset [qx_q]$.
Then $x_p=x_q$ and the triangle $[pqx_p]$ can be filled by a flat geodesic triangle.
\end{thm}

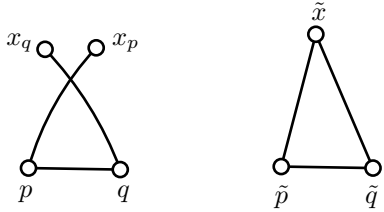
\begin{center}
\begin{lpic}[t(-0 mm),b(-0 mm),r(0 mm),l(0 mm)]{pics/key-lemma(1)}
\lbl[t]{5,1;$p$}
\lbl[t]{18,1;$q$}
\lbl[l]{16.5,20;$x_p$}
\lbl[r]{6,20;$x_q$}
\lbl[t]{39,1;$\~p$}
\lbl[t]{51,1;$\~q$}
\lbl[b]{44,23;$\~x$}
\end{lpic}
\end{center}

\parit{Proof.}
Set $p_-=p$ and $q_-=q$.

Choose  points $p_0,p_+\in [p_-,x_p]\cap \iota(N)$ so that the points $p_-,p_0,p_+,x_p$ appear in the same order on $[p_-,x_p]$.
Analogously, choose points $q_0,q_+\in [q_-,x_q]\cap N$ so that the points $q_-,q_0,q_+,x_p$ appear in the same order on $[q_-,x_q]$.
Denote by $\~p_-,\~p_0,\~p_+,\~q_-,\~q_0,\~q_+$ the corresponding points on the sides of triangle $[\~p\~q\~x]$;
so $\~p_-=\~p$ and $\~q_-=\~q$.
 
Applying the comparison to  $a_0=p_0, a_1=p_-,
a_2=p_+, a_3=x_p; \quad b_0=q_0, b_1=q_-, b_2=q_+, b_3=x_q$,
we get a model configuration $\hat p_0, \hat p_-,\hat p_+,\hat x_p,\hat q_0,\hat q_-,\hat q_+,\hat x_q$ in the Hilbert space $\HH$.

Note that from the comparison it follows that the quadruple $\hat p_-,\hat p_0,\hat p_+,\hat x_p$ lies on one line
and the same holds for the quadruple  $\hat q_-,\hat q_0,\hat q_+,\hat x_q$.

Since 
\[|\hat p_0-\hat q_+|_{\HH}\ge  |\~p_0-\~q_+|_M,\quad
|\hat p_0-\hat q_0|_{\HH}=  |\~p_0-\~q_0|_M,\quad
|\hat q_0-\hat q_+|_{\HH}=  |\~q_0-\~q_+|_M,\] 
we have $\measuredangle[\hat q_0\,{}^{\hat p_0}_{\hat q_+}]\ge \measuredangle[\~q_0\,{}^{\~p_0}_{\~q_+}]$.
The same way we get that $\measuredangle[\hat q_0\,{}^{\hat p_0}_{\hat q_-}]\ge \measuredangle[\~q_0\,{}^{\~p_0}_{\~q_-}]$.
Since the sum of adjacent angles is $\pi$, these two inequalities imply that 
\[\measuredangle[\hat q_0\,{}^{\hat p_0}_{\hat q_\pm}]= \measuredangle[\~q_0\,{}^{\~p_0}_{\~q_\pm}].\]
The same way we get that
\[\measuredangle[\hat p_0\,{}^{\hat q_0}_{\hat p_\pm}]= \measuredangle[\~p_0\,{}^{\~q_0}_{\~p_\pm}].\]

From the angle equalities, we get that 
\[|\hat p_--\hat q_+|_{\HH}\le  |\~p_--\~q_+|_M\eqlbl{eq:p-q+}\]
and the equality holds if the points $\hat p_-, \hat q_+$ lie in one plane and on the opposite sides from the line $\hat p_0 \hat q_0$.
By (3,3)-bipolar comparison the equality in \ref{eq:p-q+} indeed holds.

It follows that configuration $\hat p_0, \hat p_-,\hat p_+,\hat x_p,\hat q_0,\hat q_-,\hat q_+,\hat x_q$ is isometric to the configuration $\~ p_0, \~ p_-,\~ p_+,\~ x,\~ q_0,\~ q_-,\~ q_+,\~ x$; in particular, $\hat x_q=\hat x_p$.

By (3,3)-bipolar comparison $|x_p-x_q|_M\le |\hat x_q-\hat x_p|_{\HH}$;
therefore $x_p=x_q$; so we can set further $x=x_p=x_q$.

Note that we also proved that the angles at $p$ and $q$ in the triangle $[pqx]$ coincide with their model angles; that is,
\[\measuredangle[p\,{}^q_x]=\measuredangle[{\~p}\,{}^{\~q}_{\~x}],
\quad 
\measuredangle[q\,{}^p_x]=\measuredangle[{\~q}\,{}^{\~p}_{\~x}].\] 
By the lemma on flat slices (see for example \cite{L-extr}),
 there is a global isometric embedding $\iota'$ of the solid model triangle $[\~p\~q\~x]$ to $M$ which sends $[\~p\~q]$ to $[pq]$ and $[\~p\~x]$ to $[px]$.
Note that $\iota'$ has to coincide with $\iota$ on $N$.
It follows that $\iota'$ maps $[\~q\~x]$ to $[qx]$, which finishes the proof.
\qeds
Theorem~\ref{thm-33} 
and
Theorem~\ref{thm-mtw}
follow from the propositions below.

\begin{thm}{Proposition}\label{prop-33}
Let $M$ be a complete Riemannian manifold that satisfies (3,3)--bipolar comparison.
Then any point $p\in M$ admits a neighborhood $U\ni p$ such that if $U$ contains a nonempty open flat subset, then $U$ is flat.  
\end{thm}

\parit{Proof.}
Given a point $p$ consider a convex neighborhood $U\ni p$ such that injectivity radius at any point of $U$ exceeds the diameter of $U$; 
in particular any two points  $p,q\in U$ are connected by unique minimizing geodesic $[pq]$ which lies in~$U$.
Denote by $F$ an open flat subset in $U$; we can assume that $F$ is convex.

Note that by the key lemma we have the following:

\begin{clm}{Claim}
For any $x\in U$ and any $p,q\in F$ the triangle $[pqx]$ admits a geodesic isometric filling by a flat triangle.
\end{clm}
\medskip

Indeed, set $x_p=x$.
Consider a plane triangle $[\~p\~q\~x]$ that has the same angle at $\~p$ and the same adjacent sides as the triangle $[pqx]$.
Since $F$ is flat and convex there is a flat open geodesic surface $\Sigma$  containing $[pq]$ and a part of $[px]$ near $p$.
Choose a direction at $q$ that runs in $\Sigma$ at the angle  $\measuredangle[{\~q}\,{}^{\~p}_{\~x}]$ to $[qp]$.
Consider the geodesic  in this direction of the length
$|\~q\~x|$.
Since diameter of $U$ exceeds the injectivity radius at $q$, this geodesic is minimizing.
It remains to apply the key lemma.

From the claim, 
it follows that the sectional curvature $\sigma_x(X,Y)$ vanishes 
for any point $x\in U$ and any two velocity vectors $X,Y\in\T_x$ of minimizing geodesics from $x$ to $F$.   
Since the set of such sectional directions is open, curvature vanish at $x$; hence the result.
\qeds

\begin{thm}{Proposition}\label{prop-mtw}
Let $M$ be a complete Riemannian manifold that satisfies MTW$^{\not\perp}$.
Then any point $p\in M$ admits a neighborhood $U\ni p$ such that if $U$ contains a nonempty open flat subset, then $U$ is flat.  
\end{thm}
\parit{Proof.}
For a given $p\in M$ let us take a neighborhood $U\ni p$ as in the proof
of the previous proposition. 
The same proof as 
(Thm~1.2  
\cite{LPZ} )
shows
 that 
$U$ satisfies 
  $(4,1)$-bipolar comparison
(CTIL condition is not
necessary,
 because we stay away from 
cut-locus). Again, same proof as 
(the  Thm~1.2  
\cite{LPZ}) shows that  inside this neighborhood $(4,1)$-bipolar comparison is
equivalent to $(4,4)$-bipolar comparison.
Further note that
 $(4,4)$-bipolar comparison implies $(3,3)$-bipolar comparison.
Now we can follow the same lines as in the proof 
of Proposition~\ref{prop-33}, because  $(3,3)$-bipolar comparison
is used only locally in the proof.
\qeds

\Addresses

\end{document}